# CORRELATED SAMPLES WITH FIXED AND NONNORMAL LATENT VARIABLES[1]

By Savas Papadopoulos and Yasuo Amemiya

*Bank of Greece and IBM T. J. Watson Research Center*

A general structural equation model is fitted on a panel data set that consists of $I$ correlated samples. The correlated samples could be data from correlated populations or correlated observations from occasions of panel data. We consider cases in which the full pseudo-normal likelihood cannot be used, for example, in highly unbalanced data where the participating individuals do not appear in consecutive years. The model is estimated by a partial likelihood that would be the full and correct likelihood for independent and normal samples. It is proved that the asymptotic standard errors (a.s.e.'s) for the most important parameters and an overall-fit measure are the same as the corresponding ones derived under the standard assumptions of normality and independence for all the observations. These results are very important since they allow us to apply classical statistical methods for inference, which use only first- and second-order moments, to correlated and nonnormal data. Via a simulation study we show that the a.s.e.'s based on the first two moments have negligible bias and provide less variability than the a.s.e.'s computed by an alternative robust estimator that utilizes up to fourth moments. Our methodology and results are applied to real panel data, and it is shown that the correlated samples cannot be formulated and analyzed as independent samples. We also provide robust a.s.e.'s for the remaining parameters. Additionally, we show in the simulation that the efficiency loss for not considering the correlation over the samples is small and negligible in the cases with random and fixed variables.

Received February 2003; revised January 2005.

[1]Supported by the Department for the Supervision of Credit System and Financial Institutions of the Bank of Greece, and the Departments of Statistics of Iowa State University and Rice University.

*AMS 2000 subject classifications.* Primary 62F35, 62E20, 62H99; secondary 62H25.

*Key words and phrases.* Structural equation modeling (SEM), asymptotic robustness, multivariate analysis, panel or longitudinal data, repeated measures, nonlinear parameters, capital, credit and market risks, risk weighted assets.

**1. Introduction.** Latent variable analysis has been used widely in the social and behavioral sciences as well as in economics, and its use in medical and business applications is becoming popular. Path analysis, confirmatory factor analysis and latent variable models are the most popular psychometric models, and are all special cases of structural equation modeling (SEM). Additionally, in econometrics special cases of structural equation modeling are simultaneous equations, errors-in-variables models and dynamic panel data with random effects. In latent variable models, underlying subject-matter concepts are represented by unobservable latent variables, and their relationships with each other and with the observed variables are specified. The models that express observed variables as a linear function of latent variables are extensively used, because of their simple interpretation and the existence of computer packages such as EQS [9], LISREL [18] and PROC CALIS (SAS Institute [27]). The standard procedures in the existing computer packages assume that all the variables are normally distributed. The normality and linearity assumptions make the analysis and the interpretation simple, but their applicability in practice is often questionable. In fact, it is rather common in many applications to use the normality-based standard errors and model-fit test procedures when observed variables are highly discrete, bounded, skewed or generally nonnormal. Thus, it is of practical and theoretical interest to examine the extent of the validity of the normality-based inference procedures for nonnormal data and to explore possible ways to parameterize and formulate a model to attain wide applicability. In the structural equation analysis literature, this type of research is often referred to as asymptotic robustness study. Most existing results on this topic have been for a single sample from one population. This paper addresses the problem for multiple samples or multiple populations, and provides a unified and comprehensive treatment of the so-called asymptotic robustness. The emphasis here is the suggestion that proper parameterization and modeling lead to practical usefulness and to a meaningful interpretation. It is the first study that shows robust asymptotic standard errors (a.s.e.'s) and overall-fit measures for correlated samples with fixed factors for models with latent variables. Novel formulas are provided for the computation of the a.s.e.'s for the means and variances of the fixed correlated factors. Also, in the case of random correlated factors we prove that the a.s.e.'s of the means for the factors are robust. The superiority of the suggested a.s.e.'s to the existing robust a.s.e.'s that involve the computation of third and fourth moments is shown numerically. In a simulation study, the proposed a.s.e.'s are shown to have less variability than the robust a.s.e.'s computed by the so-called sandwich estimator. Also, the simulation studies were conducted to verify the theoretical results, assess the use of asymptotic results in finite samples, show the robustness of the power for tests and demonstrate the efficiency of the method relative to the full-likelihood estimation method that includes all the covariances of the variables over populations. The proposed method can be applied to all correlated data that can be grouped as a few correlated samples. In these correlated samples the observations are independent; for example, in panel data the correlated samples could be the occasions. The proposed methodology models variables within the samples and it can ignore the modeling of the variables between the correlated samples when it is impossible, for example, in highly unbalanced panel data in which the participating individuals do not appear in consecutive years. An application with real panel data from the Greek banking sector illustrates the importance of the proposed methodology and the derived theoretical results. In this example, it is shown that the correlated samples cannot be formulated and analyzed as independent samples.

A general latent variable model for a multivariate observation vector $\nu_j^{(i)}$ with dimension $p^{(i)} \times 1$ that is an extension of the models considered by Anderson [3, 4], Browne and Shapiro [14] and Satorra [28, 29, 30, 31, 32, 33] is

(1)
$$\nu_j^{(i)} = \beta^{(i)} + B^{(i)}\xi_j^{(i)},$$

$$\text{with } \xi_j^{(i)} = \begin{pmatrix} \zeta_j^{(i)} \\ \varepsilon_j^{(i)} \end{pmatrix} \text{ and } i = 1, \ldots, I; j = 1, \ldots, n^{(i)},$$

under the following set of assumptions. The model is extended with fixed and correlated-over-populations latent variables.

ASSUMPTION 1.

(i) There are two cases:

*Case* A: The variable $\zeta_j^{(i)}$ is (a) random with mean vector $\mu_{\zeta^{(i)}}$ and covariance matrix $\Sigma_{\zeta^{(i)}}$, (b) correlated over $i$ (i.e., the measurements of the $j$th individual of the $i_1$th population are correlated with the corresponding measurements of the $j$th individual of the $i_2$th population, for $j \leq \min\{n^{(1)}, n^{(2)}\}$) and (c) independent over $j$ (for each population the measurements of the observed individuals are independent).

*Case* B: The variable $\zeta_j^{(i)}$ is (a) fixed with limiting mean vector $\mu_{\zeta^{(i)}} = \lim_{n^{(i)} \to \infty} \bar{\zeta}^{(i)}$ and limiting covariance matrix $\Sigma_{\zeta^{(i)}} = \lim_{n^{(i)} \to \infty} \mathbf{S}_{\zeta^{(i)}}$ and (b) correlated over $i$ [see comments in case A(b)].

(ii) There exists $\varepsilon_j^{(i)} = (\varepsilon_{0j}^{(i)\prime}, \varepsilon_{1j}^{(i)\prime}, \ldots, \varepsilon_{L^{(i)}j}^{(i)\prime})'$, where (a) $\varepsilon_{0j}^{(i)} \sim N(0, \Sigma_{\varepsilon_0^{(i)}})$, (b) $\varepsilon_{\ell j}^{(i)}$ ($\ell = 1, \ldots, L^{(i)}$) are independent over $i, \ell$ and $j$ with mean 0 and covariance matrix $\Sigma_{\varepsilon_\ell^{(i)}}$ and (c) $\zeta_j^{(i)}$ are independent with $\varepsilon_{\ell j}^{(i)}$ ($\ell = 0, 1, \ldots, L^{(i)}$) over $i$ and $j$.

(iii) The intercepts $\beta^{(i)}$, the coefficients $B^{(i)}$ and the variance matrices of the normally distributed errors $\Sigma_{\varepsilon_0^{(i)}}$ can be restricted. Thus, they are assumed to be functions of a vector $\tau$.

(iv) The mean vectors $\mu_{\zeta^{(i)}}$, the variance matrices $\Sigma_{\zeta^{(i)}}$ of the correlated factors and the variance matrices of the nonnormal vectors $\Sigma_{\varepsilon_\ell^{(i)}}$ ($\ell = 1, \ldots, L^{(i)}$) are assumed to be unrestricted.

A common approach to verifying the identification and fitting the model is to assume hypothetically that all $\xi_j^{(i)}$'s are normally distributed and to concentrate on the first two moments of the observed vector $\nu_j^{(i)}$. The issue for the so-called asymptotic robustness study is to assess the validity of such procedures based on the assumed normality, in terms of inference for unknown parameters, for a wide class of distributional assumptions on $\xi_j^{(i)}$. It turns out that the type of parameterization used in the model, restricting the coefficient $B^{(i)}(\tau)$ but keeping the variances $\Sigma_{\varepsilon_\ell^{(i)}}$ of the nonnormal latent variables $\varepsilon_{\ell j}^{(i)}$ unrestricted, plays a key role in the study.

The model, the notation and the assumptions are explained by the following example.

EXAMPLE 1.  A two-population ($I = 2$) recursive system of simultaneous equations with errors in the explanatory variables is considered. The model is shown in (2). The system in (2) can be written in the matrix form $\nu_j^{(i)} = \alpha^{(i)} + \Gamma^{(i)}\nu_j^{(i)} + \Delta^{(i)}\zeta_j^{(i)} + \mathbf{e}_j^{(i)}$, which has the form of model (1) with $\beta^{(i)} = (\mathbf{I}^{(i)} - \Gamma^{(i)})^{-1}\alpha^{(i)}$, $B^{(i)} = (\mathbf{I}^{(i)} - \Gamma^{(i)})^{-1}[\Delta^{(i)}, \mathbf{I}^{(i)}]$ and $\varepsilon_j^{(i)} = \mathbf{e}_j^{(i)}$. The model is also a special case of the LISREL model with no latent variables in the dependent variables $\mathbf{y}^{(i)}$, that is, $\mathbf{y}^{(i)} = \eta^{(i)}$, in the LISREL notation. The latent variables $\zeta_j^{(1)}$ and $\zeta_j^{(2)}$ are correlated for each $j = 1, \ldots, 500$, with correlation 0.4. That is, the measurements of each individual from the second population are correlated with the measurements of one individual from the first population. The first population also has 500 individuals that are independent from all the individuals of the second population. Note that the number of observed variables is different for the two populations. Four measurements, $x_j^{(1)}, y_{1j}^{(1)}, y_{2j}^{(1)}$ and $y_{3j}^{(1)}$, are taken from the first population ($p^{(1)} = 4$) and three measurements, $x_j^{(2)}, y_{1j}^{(2)}$ and $y_{2j}^{(2)}$, are taken from the second ($p^{(2)} = 3$). For $j = 1, \ldots, n^{(i)}$, with $n^{(1)} = 1000$ and $n^{(2)} = 500$,

$$x_j^{(1)} = \zeta_j^{(1)} + e_{0j}^{(1)}, \qquad\qquad x_j^{(2)} = \zeta_j^{(2)} + e_{0j}^{(2)},$$
$$y_{1j}^{(1)} = \beta_1 + \delta_1\zeta_j^{(1)} + e_{1j}^{(1)}, \qquad\qquad y_{1j}^{(2)} = \beta_1 + \delta_1\zeta_j^{(2)} + e_{1j}^{(2)},$$

(2)
$$y_{2j}^{(1)} = \beta_2 + \gamma_1 y_{1j}^{(1)} + \delta_2 \zeta_j^{(1)} + e_{2j}^{(1)}, \qquad y_{2j}^{(2)} = \beta_2 + \gamma_1 y_{1j}^{(2)} + \delta_2 \zeta_j^{(2)} + e_{2j}^{(2)},$$

$$y_{3j}^{(1)} = \beta_3 + \gamma_2 y_{2j}^{(1)} + e_{3j}^{(1)}.$$

The parameters $\beta_1, \beta_2, \gamma_1, \delta_1$ and $\delta_2$ do not depend on $i$. That is, they are common for the two populations. These parameters belong to the vector $\tau$. The variables $\zeta_j^{(1)}$ and $\zeta_j^{(2)}$ can be fixed or nonnormal according to cases A and B of Assumption 1. If all the errors are normal in accordance with the notation of Assumption 1, we have $\varepsilon_{0j}^{(i)} = \mathbf{e}_j^{(i)}$, while if $\mathbf{e}_{0j}^{(i)}$ is normal and all the other errors are nonnormal, then $\varepsilon_{0j}^{(i)} = e_{0j}^{(i)}$ and $\varepsilon_{\ell j}^{(i)} = e_{\ell j}^{(i)}$ for $i = 1, 2$, $j = 1, \ldots, n^{(i)}$ and $\ell = 1, \ldots, L^{(i)}$ with $L^{(1)} = 3$ and $L^{(2)} = 2$. According to Assumption 1, only the variances of the normal errors can be restricted to be the same over populations and these variances belong to the vector $\tau$.

Further discussion about the model in (1) is given in Section 2. The model in (2) of Example 1 is simulated in Section 4 and used as an example to explain the theory in this paper.

Latent variable analysis of multiple populations was discussed by Jöreskog [17], Lee and Tsui [20], Muthén [23] and Satorra [29, 30]. The so-called asymptotic robustness of normal-based methods for latent variable analysis has been extensively studied in the last 15 years. For exploratory (unrestricted) factor analysis, Amemiya, Fuller and Pantula [2] proved that the limiting distribution of some estimators is the same for fixed, nonnormal and normal factors under the assumption that the errors are normally distributed. Browne [12] showed that the above results hold for a more general class of latent variable models assuming finite eighth moments for the factors and normal errors. Anderson and Amemiya [5], and Amemiya and Anderson [1] extended the above results to confirmatory factor analysis and nonnormal errors; they assumed finite second moments for the factors and errors. Browne and Shapiro [14] introduced a general linear model and used an approach based on the finite fourth moments that differs from that of Anderson and Amemiya. Considering the model of Browne and Shapiro, Anderson [3, 4] included nonstochastic latent variables and assumed only finite second moments for the nonnormal latent variables. Latent variable models with mean and covariance structures were studied by Browne [13] and Satorra [28]. Satorra [29, 30, 31, 32, 33] first considered asymptotic robustness for linear latent models in multisample analysis of augmented-moment structures. Additional studies on the asymptotic robustness of latent variable analysis were conducted by Shapiro [37], Mooijaart and Bentler [22] and Satorra and Bentler [35].

For the one-sample problem, asymptotic distribution-free (ADF) methods for latent variable analysis were proposed to deal with nonnormal data (see, e.g., [8, 11, 23]). The ADF methods turned out to be problematic in practice, since the fourth-order sample moments are very variable (see, e.g., [15, 24]). In this paper mean and covariance structures are considered for a general multipopulation model that contains fixed, normal and nonnormal variables; some of the nonnormal variables are allowed to be correlated over populations. We use the approach of Anderson and Amemiya [5] to show that the normal-based methods are applicable for nonnormal and nonrandom data assuming finite second-order moments. We also use extensively theory and notation from matrix analysis (see, e.g., [16, 21]).

Section 2 explains the suggested parameterization and the estimation procedure. The theoretical results are derived and discussed in Section 3. Section 4 reports results from simulation studies and that the proposed asymptotic standard errors seem to be numerically more efficient than those derived by the sandwich estimator. Our methodology and the theoretical results are applied and explained in Section 5 by the fit of an econometric model with latent economic factors to real data.

## 2. Model, parameterization and procedure.

In this paper we study the model (1) introduced in Section 1. We consider $I$ populations and we assume that $n^{(i)}$ individuals are sampled from the $i$th population, $i = 1, \ldots, I$, and that $p^{(i)}$ measurements are taken from each sampled individual in the $i$th population. Denote the multisample data set by $\nu_j^{(i)}, i = 1, \ldots, I, j = 1, \ldots, n^{(i)}$, where $\nu_j^{(i)}$ is the $p^{(i)} \times 1$ measurement vector from the $j$th individual in the $i$th population. We consider a very general latent variable model that includes models widely used in single population cases and covers a large class of distributional situations in one form. To cover various distributional settings, it is convenient to assume that the observed vector $\nu_j^{(i)}$ can be written as a linear combination of $L^{(i)} + 2$ independent latent vectors and that the latent vectors can be divided into three groups: (1) a fixed or nonnormal vector that is correlated over populations $\zeta_j^{(i)}$, (2) a random vector $\varepsilon_{0j}^{(i)}$ assumed to be normally distributed and (3) $L^{(i)}$ nonnormal vectors $\varepsilon_{\ell j}^{(i)}$ $(\ell = 1, \ldots, L^{(i)})$. Note that the sample size $n^{(i)}$, the number of measured variables $p^{(i)}$ and the number of latent vectors $L^{(i)}$ generally differ over populations (dependent on $i$). The generality of this model allows us to deal with cases where slightly different variables are measured from different populations with possibly different structures.

All normally distributed latent variables are included in $\varepsilon_{0j}^{(i)}$ and their distribution may possibly be related through $\tau$ over populations $i = 1, \ldots, I$. Other unspecified or nonnormal random latent variables are divided into independent parts $\ell = 1, \ldots, L^{(i)}$ with unrestricted covariance matrices. Case

A of Assumption 1 with fixed $\zeta_j^{(i)}$ can represent a situation where the interest is in the model fitting and estimation only for a given set of individuals and not for the populations. In addition, the fixed $\zeta_j^{(i)}$ can be used in an analysis conducted conditionally on a given set of $\zeta_j^{(i)}$ values. Such a conditional analysis may be appropriate when the individuals $j = 1, \ldots, n^{(i)}$ do not form a random sample from the $i$th population and/or when a component of $\nu_j^{(i)}$ represents some dependency over $I$ populations. For example, the $I$ populations may actually correspond to a single population at $I$ different time points. With $\zeta_j^{(i)}$ being latent and fixed, the limits of the unobservable sample mean, $\mu_{\zeta^{(i)}}$, and of the sample covariance matrix, $\Sigma_{\zeta^{(i)}}$, are assumed to be unknown and unrestricted. All $\beta^{(i)}(\tau)$ and $B^{(i)}(\tau)$ are expressed in terms of $\tau$, which represents known or restricted elements and allows functional relationships over $I$ populations. Even though $\tau$ also appears in $\Sigma_{\varepsilon_0^{(i)}}(\tau)$, the elements of $\tau$ are usually divided into two groups: one for $\Sigma_{\varepsilon_0^{(i)}}(\tau)$, and another for $\beta^{(i)}(\tau)$ and $B^{(i)}(\tau)$.

Assumption 1(iii) and (iv) provide a particular identifiable parameterization for the model in (1). For the single population case with $I = 1$, various equivalent parameterizations have been used in practice. Some place restrictions on covariance matrices (e.g., by standardizing latent variables) and leave the coefficients unrestricted. The parameterization that leaves the covariance matrices (and possibly some mean vectors) of latent variables unrestricted and that places identification restrictions only on the coefficients and intercepts is referred to as the errors-in-variables parameterization. For the single population case, a parameterization with restricted covariance matrices generally has an equivalent errors-in-variables parameterization, and the two parameterizations with one-to-one correspondence lead to an equivalent interpretation. The one-sample asymptotic robustness results have shown that the asymptotic standard errors for the parameters in the errors-in-variables formulation computed under the normality assumption are valid for nonnormal data, but that the same does not hold under parameterization with restricted covariance matrices. For the multisample, the model in (1), we will show that the errors-in-variables type parameterization given in Assumption 1 provides asymptotic robustness. However, for the multisample case there are other reasons to consider the parameterization specified in Assumption 1(iii) and (iv). As mentioned earlier, a multipopulation study is conducted because the populations are thought to be different, but certain aspects of the structure generating data are believed to be common over populations. Suppose that the same or similar measurements are taken from different populations. For example, a similar set of psychological tests may be given to a number of different groups, for example, two gender groups, groups with different occupations or educational backgrounds, groups in varying socioeconomic or cultural environments, or different time points in the growth of a group. The subject matter or scientific interest exists in making inferences about some general assertion that holds commonly for various populations. Such interest is usually expressed as relationships among latent (and observed) variables that hold regardless of the location and variability of the variables. Then a relevant analysis is to estimate and test the relationships, and to explore the range of populations for which the relationships hold. The parameterization in Assumption [1](iii) and (iv) with unrestricted $\Sigma_{\varepsilon_\ell^{(i)}}$ and generally structured $B^{(i)}(\tau)$ corresponds very well with the scientific interest of the study, and allows an interpretation consistent with the practical meaning of the problem. Note that $\Sigma_{\varepsilon_\ell^{(i)}}, i = 1, \ldots, I, \ell = 1, \ldots, L^{(i)}$, are unrestricted covariance matrices and do not have any relationships over $i$ or $\ell$, and that $\beta^{(i)}(\tau)$ and $B^{(i)}(\tau)$ can have known elements and elements with relationships over $i$ and $\ell$. On the other hand, the covariance matrix $\Sigma_{\varepsilon_0^{(i)}}$ of the normal latent vector $\varepsilon_{0j}^{(i)}$ can have restrictions or equality over populations through $\tau$. This gives the generality of the model in [(1)](1) with only one normal latent vector, because a block diagonal $\Sigma_{\varepsilon_0^{(i)}}$ corresponds to a number of independent subvectors in the normal $\varepsilon_{0j}^{(i)}$. In addition, the possibility of restrictions on $\Sigma_{\varepsilon_0^{(i)}}$ over populations can also be important in applications. For example, if the same measurement instruments are applied to different samples, then the variances of pure measurement errors may be assumed to be the same over the samples. However, the normality assumption for pure measurement errors is reasonable in most situations, and such errors can be included in $\varepsilon_{0j}^{(i)}$. Assumption [1](iv) and (v) do not rule out latent variable variances and covariances with restrictions across populations, but do require the latent variables with restricted variances to be normally distributed. This requirement is not very restrictive in most applications, as discussed above, but it is needed to obtain the asymptotic robustness results given in the next section. The general form of $\beta^{(i)}(\tau)$ and inclusion of the fixed latent vector allow virtually any structure for the means of the observed $\nu_j^{(i)}$. Hence, the errors-in-variables type parameterization in Assumption [1](iii) can solve the identification problem, provide a general and convenient way to represent the subject-matter theory and concepts, and produce asymptotic robustness results presented in the next section.

For the multisample data $\nu_j^{(i)}$ in [(1)](1), let $\bar{\nu}^{(i)}$ and $\mathbf{S}_\nu^{(i)}$ be the sample mean vector and sample covariance matrix (unbiased) for the $i$th population, $i = 1, \ldots, I$. It is assumed that the sample covariance matrices $\mathbf{S}_\nu^{(i)}$ are nonsingular with probability 1. Define

$$
(3) \qquad \mathbf{c}^{(i)} = \begin{pmatrix} \bar{\nu}^{(i)} \\ \mathrm{vec}(\mathbf{S}_\nu^{(i)}) \end{pmatrix}, \qquad \mathbf{c} = \begin{pmatrix} \mathbf{c}^{(1)} \\ \vdots \\ \mathbf{c}^{(I)} \end{pmatrix}.
$$

We consider model fitting and estimation based only on $\mathbf{c}$, because such procedures are simple and have some useful properties. Also note that Assumption 1 does not specify a particular distributional form of observations beyond the first two moments and specifies no particular correspondence or relationship between samples. Let $\theta$ be a $d_\theta \times 1$ vector containing all unknown parameters in $E(\mathbf{c}) = \gamma(\theta)$ under the model in (1) and Assumption 1, and let $\theta = (\tau', \upsilon')'$, where $\tau$ and $\upsilon$ contain the parameters mentioned in Assumption 1(iv) and (v), respectively. That is, $\tau$ contains parameters that can be restricted, while $\upsilon$ contains the parameters that cannot be restricted over populations. Under the model in (1) and Assumption 1, we compute the expected means

$$
\mu_\nu^{(i)}(\theta) = E(\bar{\nu}^{(i)}) \quad \text{and} \quad \Sigma_\nu^{(i)}(\theta) = E(\mathbf{S}_\nu^{(i)}).
$$

For the estimation of $\theta$, we consider an estimator $\hat{\theta}$ obtained by minimizing over the parameter space

$$
(4) \qquad \begin{aligned}
Q(\theta) = \sum_{i=1}^I n^{(i)} \{ &\mathrm{tr}[\mathbf{S}_\nu^{(i)} \Sigma_\nu^{(i)-1}(\theta)] - \log|\mathbf{S}_\nu^{(i)} \Sigma_\nu^{(i)-1}(\theta)| - p^{(i)} \\
&+ [\bar{\nu}^{(i)} - \mu_\nu^{(i)}(\theta)]' \Sigma_\nu^{(i)-1}(\theta)[\bar{\nu}^{(i)} - \mu_\nu^{(i)}(\theta)] \}.
\end{aligned}
$$

The obtained estimator $\hat{\theta}$ is a slight modification of the normal maximum likelihood estimator (MLE). The exact normal MLE can be obtained if $[(n^{(i)} - 1)/n^{(i)}]\mathbf{S}_\nu^{(i)}$ is used in place of $\mathbf{S}_\nu^{(i)}$. Asymptotic results are equivalent for the two estimators. We consider $\hat{\theta}$ because it can be computed with existing computer packages. The form of $Q(\theta)$ corresponds to the so-called mean and covariance structure analysis, but the existing covariance structure computer packages without mean structure can be used to carry out the minimization of $Q(\theta)$ using a certain technique (see, e.g., the EQS and LISREL manuals). Note that other estimation techniques that are asymptotically equivalent to MLE can be used, such as minimum distance, which is actually a generalization of the generalized method of moments. In the next section, asymptotic distribution results for $\hat{\theta}$ are derived for a broad range of situations.

**3. Theoretical results.** The main results of this paper are presented in Theorem 1. We now define a set of assumptions for the model in (1) that assumes normal and independent variables over populations under the same parameterization as in Assumption 1.

ASSUMPTION 1B.

(i) For all $i$ and $j$ $(i = 1, \ldots, I; j = 1, \ldots, n^{(i)})$ $\zeta_j^{(i)} \sim N(\mu_{\zeta^{(i)}}, \Sigma_{\zeta^{(i)}})$ and are independent.

(ii) For all $\ell = 0, 1, \ldots, L^{(i)}, \varepsilon_\ell^{(i)} \sim N(0, \Sigma_{\varepsilon_\ell^{(i)}})$.

(iii) The matrices $\beta^{(i)}, B^{(i)}$ and $\Sigma_{\varepsilon_0^{(i)}}$ can be restricted and are assumed to be functions of a vector $\tau$.

(iv) The matrices $\mu_{\zeta^{(i)}}, \Sigma_{\zeta^{(i)}}$ and $\Sigma_{\varepsilon_\ell^{(i)}}, \ell = 1, \ldots, L^{(i)}$, are assumed to be unrestricted.

Theorem 1 shows similarities and differences of the limiting results for the two different sets of Assumptions 1 and 1B.

THEOREM 1. *Assume that the model in* (1) *holds under Assumption* 1. *In addition we make the following assumptions:*

ASSUMPTION 2. There exists $\lim_{n_m \to \infty}(n^{(i)}/n) = r^{(i)}$, where $n_m = \min\{n^{(1)}, \ldots, n^{(I)}\}$ and $n = \sum_{i=1}^I n^{(i)}$.

ASSUMPTION 3. $(\forall \varepsilon > 0)(\exists \delta > 0) \ni |\gamma(\theta) - \gamma(\theta_0)| < \delta \Rightarrow \|\theta - \theta_0\| < \varepsilon$, where $\|\mathbf{x}\| = \sqrt{\mathbf{x}'\mathbf{x}}$ and $\theta_0$ is the limiting true value of $\theta$.

ASSUMPTION 4. For all $i = 1, \ldots, I$, $\beta^{(i)}(\tau)$, $B^{(i)}(\tau)$ and $\Sigma_{\varepsilon_0^{(i)}}(\tau)$ are twice continuously differentiable in the parameter space of $\tau$. The columns of the matrix $\partial\gamma(\theta_0)/\partial\tau'$ are linearly independent.

THEOREM 1 (cont.).

(i) *Then*

$$\mathbf{V}_{\mathrm{G}}^{(\tau)} = \mathbf{V}_{\mathrm{NI}}^{(\tau)},$$

*where* $\mathbf{V}_{\mathrm{G}}^{(\tau)}$ *and* $\mathbf{V}_{\mathrm{NI}}^{(\tau)}$ *are the asymptotic covariance matrices of* $\hat{\tau}$ *under the general Assumption* 1 *and under the standard Assumption* 1B, *respectively (the initials NI stand for normality and independence over populations and G stands for the general set of Assumptions* 1). *The matrix* $\mathbf{V}_{\mathrm{G}}^{(\tau)}$ *is the part of the matrix* $\mathbf{V}_{\mathrm{G}}^{(\theta)}$ *that is the asymptotic covariance matrix for the estimated vector* $\theta$.

(ii) *For the asymptotic covariance matrices for the mean vectors* $\widehat{\mu}_{\zeta^{(i)}}$, (1) *in case* A *of Assumption* 1 *with fixed* $\zeta_j^{(i)}$,

$$\mathbf{V}_{\mathrm{G}}^{(\mu_{\zeta^{(i)}})} = \mathbf{V}_{\mathrm{NI}}^{(\mu_{\zeta^{(i)}})} - \Sigma_{\zeta^{(i)}} \tag{5}$$

*holds, and* (2) *in case* B *of Assumption* 1 *with random* $\zeta_j^{(i)}$,

$$\mathbf{V}_{\mathrm{G}}^{(\mu_{\zeta^{(i)}})} = \mathbf{V}_{\mathrm{NI}}^{(\mu_{\zeta^{(i)}})} \tag{6}$$

*holds.*

(iii) *For the asymptotic covariance matrices for* $\mathrm{vec}(\widehat{\Sigma}_{\zeta^{(i)}})$, (1) *in case* A *of Assumption* 1 *with fixed* $\zeta_j^{(i)}$,

$$\mathbf{V}_{\mathrm{G}}^{(\mathrm{vec}(\Sigma_{\zeta^{(i)}}))} = \mathbf{V}_{\mathrm{NI}}^{(\mathrm{vec}(\Sigma_{\zeta^{(i)}}))} - \frac{2}{n^{(i)}}(\Sigma_{\zeta^{(i)}} \otimes \Sigma_{\zeta^{(i)}}) \tag{7}$$

*holds, and* (2) *in case* B *of Assumption* 1 *with random* $\zeta_j^{(i)}$ *and assuming that* $\zeta_j^{(i)}$ *have finite fourth moments,*

$$\mathbf{V}_{\mathrm{G}}^{(\mathrm{vec}(\Sigma_{\zeta^{(i)}}))} = \mathbf{V}_{\mathrm{NI}}^{(\mathrm{vec}(\Sigma_{\zeta^{(i)}}))} + \frac{1}{n^{(i)}}\mathrm{Var}[\mathrm{vec}(\zeta^{(i)}\zeta^{(i)\prime})] - \frac{2}{n^{(i)}}(\Sigma_{\zeta^{(i)}} \otimes \Sigma_{\zeta^{(i)}}) \tag{8}$$

*holds.*

(iv) *The function* $Q(\theta)$, *defined in* (4), *evaluated on its minimum* $\widehat{\theta}$ *converges to a chi-square distribution,* $Q(\widehat{\theta}) \xrightarrow{d} \chi_q^2$, *with* $q = \sum_{i=1}^{I}[p^{(i)} + p^{(i)}(p^{(i)} + 1)/2] - d_\theta$.

PROOF OF THEOREM 1. For the proof we need the following three lemmas.

LEMMA 1. *Assume that the model in* (1) *holds. If Assumptions* 1, 2 *and* 3 *hold, then as* $n_m \to \infty$,

$$\widehat{\theta} \xrightarrow{p} \theta_0. \tag{9}$$

PROOF. From Assumption 1 and the law of large numbers, $\mathbf{c} \xrightarrow{p} \gamma(\theta_0)$, which implies $Q(\theta_0) \xrightarrow{p} 0$. Since $Q(\theta) > 0 \; \forall \theta$ and $\widehat{\theta}$ minimizes $Q$, we have $Q(\widehat{\theta}) \xrightarrow{p} 0$. From the last result and Assumption 2 we get $\gamma(\widehat{\theta}) \xrightarrow{p} \gamma(\theta_0)$, and (9) holds from Assumption 3. $\square$

LEMMA 2. *Let* $\bar{\theta}_n = (\tau_0', v_n')'$, *where* $\tau_0$ *is the true value of* $\tau$ *and* $v_n$ *contains the vectors* $\bar{\zeta}^{(i)}$, $\mathrm{vec}(\mathbf{S}_{\zeta^{(i)}})$ *and* $\mathrm{vec}(\mathbf{S}_{\varepsilon_\ell^{(i)}})$, $\ell = 1, \ldots, L^{(i)}$, *for all* $i = 1, \ldots, I$.

(i) *Then, under the model and the assumptions considered in Lemma* 1, *and under Assumption* 4,

$$\sqrt{n}(\hat{\theta} - \bar{\theta}_n) = \mathbf{A}_0 \sqrt{n} [\mathbf{c} - \gamma(\bar{\theta}_n)] + o_p(1), \tag{10}$$

*where* $\mathbf{A}_0$ *is free of* $n^{(i)}$ *and*

$$\mathbf{A}_0 = (\mathbf{J}_0' \Omega_0^{-1} \mathbf{J}_0)^{-1} \mathbf{J}_0' \Omega_0^{-1}, \tag{11}$$

*where* $\mathbf{J}_0 = \mathbf{J}(\gamma(\theta_0))$ *is the Jacobian of* $\gamma(\theta)$ *evaluated at* $\theta_0$, $\Omega_0^{-1} = \Omega^{-1}(\theta_0) = [r^{(1)} \Omega^{(1)-1}(\theta_0)] \oplus \cdots \oplus [r^{(I)} \Omega^{(I)-1}(\theta_0)]$ *and* $\Omega^{(i)-1}(\theta) = \Sigma^{(i)-1}(\theta) \oplus \{\frac{1}{2} \times [\Sigma^{(i)-1}(\theta) \otimes \Sigma^{(i)-1}(\theta)]\}$.

*Recall that the ratios* $r^{(i)}$ *were defined in Assumption* 2 *and* $\mathbf{c}$ *was defined in* (3). *The symbol* $\oplus$ *is the direct sum for matrices.*

(ii) *Also,*

$$Q(\hat{\theta}) = n[\mathbf{c} - \gamma(\bar{\theta}_n)]' \mathbf{M}_0 [\mathbf{c} - \gamma(\bar{\theta}_n)] + o_p(1) \tag{12}$$

*with* $\mathbf{M}_0 = \Omega_0^{-1}(\mathbf{I} - \mathbf{A}_0)$.

PROOF. (i) From Taylor's expansion and Lemma 1 it turns out that there exists $\theta^*$ on the line segment between $\hat{\theta}$ and $\bar{\theta}_n$ such that

$$\mathbf{J}'[Q(\hat{\theta})] = \mathbf{J}'[Q(\bar{\theta}_n)] + \mathbf{H}[Q(\theta^*)](\hat{\theta} - \bar{\theta}_n), \tag{13}$$

where $\mathbf{J}$ and $\mathbf{H}$ are the Jacobian and Hessian matrices, respectively. Now for the Jacobian and Hessian matrices we proved that

$$\mathbf{J}'[Q(\bar{\theta}_n)] = -2\mathbf{J}_0' \Omega_0^{-1} [\mathbf{c} - \gamma(\bar{\theta}_n)] + o_p(n_m^{-1/2}), \tag{14}$$

$$\mathbf{H}[Q(\theta^*)] \xrightarrow{p} 2\mathbf{J}_0' \Omega_0^{-1} \mathbf{J}_0. \tag{15}$$

The result in (10) follows if we use (14), (15) and the fact that $\mathbf{J}[Q(\hat{\theta})] = 0$ in (13).

(ii) After doing several matrix modifications, we get the quadratic form

$$Q(\hat{\theta}) = n[\mathbf{c} - \gamma(\hat{\theta})]' \Omega^{-1}(\hat{\theta})[\mathbf{c} - \gamma(\hat{\theta})] + o_p(1). \tag{16}$$

Also, there exists $\theta^*$ on the line segment between $\hat{\theta}$ and $\bar{\theta}_n$ such that

$$\gamma(\hat{\theta}) - \gamma(\bar{\theta}_n) = \mathbf{J}[\gamma(\theta^*)](\hat{\theta} - \bar{\theta}_n). \tag{17}$$

From (17) and (10) we get that

$$\mathbf{c} - \gamma(\hat{\theta}) = [\mathbf{I} - \mathbf{J}_0 \mathbf{A}_0][\mathbf{c} - \gamma(\bar{\theta}_n)] + o_p\left(\frac{1}{\sqrt{n}}\right), \tag{18}$$

and the result follows from (16) and (18). □

LEMMA 3. (i) *For the model in* (1) *under Assumption* 1 *it holds that*

$$\mathbf{c} - \gamma(\bar{\theta}_n) = \mathbf{E}\mathbf{w}, \tag{19}$$

*where* $\mathbf{E}$ *is a constant matrix,* $\mathbf{w}$ *consists of the subvectors* $\mathbf{w}^{(i)}, i = 1, \ldots, I,$ *and* $\mathbf{w}^{(i)}$ *consists of the subvectors* $\bar{\varepsilon}^{(i)}, \mathrm{vec}(\mathbf{S}_{\varepsilon_0^{(i)}\varepsilon_0^{(i)}})$ *and* $\mathrm{vec}(\mathbf{S}_{\mathbf{x}^{(i)}\mathbf{y}^{(i)}})$ *for all* $\mathbf{x}^{(i)}$ *and* $\mathbf{y}^{(i)}$ *such that* $\mathbf{x}^{(i)} \neq \mathbf{y}^{(i)}, i = 1, \ldots, I,$ *and* $\mathbf{x}^{(i)}, \mathbf{y}^{(i)} = \zeta^{(i)}, \varepsilon_0^{(i)}, \varepsilon_1^{(i)}, \ldots,$ $\varepsilon_{L^{(i)}}^{(i)}$.

(ii) *The limiting distribution of* $\sqrt{n}\mathbf{w}$ *is the same under Assumptions* 1 *and* 1B.

PROOF. (i) We proved that the components of $\mathbf{c} - \gamma(\bar{\theta}_n)$ are written in the form

$$\bar{\nu}^{(i)} - \mu_\nu^{(i)}(\bar{\theta}_n) = B^{(i)} \begin{bmatrix} \mathbf{0} \\ \bar{\varepsilon}^{(i)} \end{bmatrix}, \tag{20}$$

$$\mathbf{S}_{\nu^{(i)}} - \Sigma_{\nu^{(i)}}(\bar{\theta}_n) = B^{(i)} \begin{bmatrix} \mathbf{0} & \mathbf{S}_{\zeta^{(i)}\varepsilon^{(i)}} \\ \mathbf{S}_{\varepsilon^{(i)}\zeta^{(i)}} & \mathbf{S}_{\varepsilon^{(i)}\varepsilon^{(i)}} - \mathbf{D}_{\varepsilon^{(i)}} \end{bmatrix} B'^{(i)}, \tag{21}$$

where $\mathbf{D}_{\varepsilon^{(i)}} = \mathbf{0} \oplus \mathbf{S}_{\varepsilon_1^{(i)}} \oplus \cdots \oplus \mathbf{S}_{\varepsilon_{L^{(i)}}^{(i)}}$. The result in (19) follows by noting in (20) and (21) that the components of $\mathbf{c} - \gamma(\bar{\theta}_n)$ are products of constant matrices (functions of $B^{(i)}$) and the subvectors of $\mathbf{w}^{(i)}$, and also using the property $\mathrm{vec}(\mathbf{ABC}) = (\mathbf{C}' \otimes \mathbf{A})\mathrm{vec}(\mathbf{B})$.

(ii) Note that the matrix $\mathbf{S}_{\varepsilon^{(i)}\varepsilon^{(i)}} - \mathbf{D}_{\varepsilon^{(i)}}$ does not depend on $\mathbf{S}_{\varepsilon_\ell^{(i)}\varepsilon_\ell^{(i)}}$ for $\ell = 1, \ldots, L^{(i)}$. Also note that within the populations for each (i) the subvectors of $\sqrt{n}\mathbf{w}^{(i)}$ are independent and their limiting distributions do not depend on the nonnormality of the latent variables and on the fixed latent variables in case A (see [4], Theorem 5.1). Now between the populations, the limiting covariance between $\mathbf{w}^{(i)}$ and $\mathbf{w}^{(m)}$ for $i \neq m$ is 0 despite the correlation of $\zeta_j^{(i)}$ and $\zeta_j^{(m)}$ for each $j$. This holds because the limiting covariance between $\sqrt{n}\,\mathrm{vec}(\mathbf{S}_{\zeta^{(i)}\varepsilon^{(i)}})$ and $\sqrt{n}\,\mathrm{vec}(\mathbf{S}_{\zeta^{(m)}\varepsilon^{(m)}})$ is 0 since the errors are assumed to be independent over populations. □

Now we return to the proof of Theorem 1. For (i) Lemmas 2(i) and 3(i) show that $\sqrt{n}(\hat{\tau} - \tau_0)$ is a linear combination of $\sqrt{n}\mathbf{w}$ and thus the result follows from Lemma 3(ii).

For cases (ii) and (iii) we use the respective equations

$$\sqrt{n}(\hat{\mu}_{\zeta^{(i)}} - \mu_{\zeta^{(i)}}^0) = \sqrt{n}(\hat{\mu}_{\zeta^{(i)}} - \bar{\zeta}^{(i)}) + \sqrt{n}(\bar{\zeta}^{(i)} - \mu_{\zeta^{(i)}}^0), \tag{22}$$

$$\sqrt{n}\mathrm{vec}(\hat{\Sigma}_{\zeta^{(i)}} - \Sigma_{\zeta^{(i)}}^0) = \sqrt{n}\mathrm{vec}(\hat{\Sigma}_{\zeta^{(i)}} - \mathbf{S}_{\zeta^{(i)}}) + \sqrt{n}\mathrm{vec}(\mathbf{S}_{\zeta^{(i)}} - \Sigma_{\zeta^{(i)}}^0), \tag{23}$$

where $\mu_{\zeta^{(i)}}^0$ and $\Sigma_{\zeta^{(i)}}^0$ are the true values of the corresponding parameters. In both (ii) and (iii), for case A with fixed factors, we need the limiting distributions of the first vectors in the second parts of (22) and (23). For case B with random factors, we need the limiting distributions of the vectors in the first parts of (22) and (23). Since the procedure is the same for (ii) and (iii), we explain the proof only for part (iii). So for case A in (23) we compute the limiting covariance matrices of all three vectors under the Assumption 1B,

$$(24) \qquad \mathbf{V}_{\mathrm{NI}}^{(\mathrm{vec}(\Sigma_{\zeta^{(i)}}))} = \mathbf{V}_2 + \frac{2}{n^{(i)}}(\Sigma_{\zeta^{(i)}} \otimes \Sigma_{\zeta^{(i)}}).$$

From Lemmas 2(i) and 3 it follows that the first vector of the second part of (23) has the same limiting distribution under Assumption 1 with fixed factors and under Assumption 1B. Thus $\mathbf{V}_2 = \mathbf{V}_{\mathrm{G}}^{(\mathrm{vec}(\Sigma_{\zeta^{(i)}}))}$ and the result follows by solving (24) for $\mathbf{V}_{\mathrm{G}}^{(\mathrm{vec}(\Sigma_{\zeta^{(i)}}))}$.

Now for case B in (iii) we compute the limiting covariance matrices under Assumption 1B and under Assumption 1, and we get, respectively,

$$(25) \qquad \mathbf{V}_{\mathrm{NI}}^{(\mathrm{vec}(\Sigma_{\zeta^{(i)}}))} = \mathbf{V}_{\mathrm{NI}}^* + \frac{2}{n^{(i)}}(\Sigma_{\zeta^{(i)}} \otimes \Sigma_{\zeta^{(i)}}),$$

$$(26) \qquad \mathbf{V}_{\mathrm{G}}^{(\mathrm{vec}(\Sigma_{\zeta^{(i)}}))} = \mathbf{V}_{\mathrm{G}}^* + \frac{1}{n^{(i)}}\mathrm{Var}[\mathrm{vec}(\zeta^{(i)}\zeta^{(i)'})].$$

Again, from Lemmas 2 and 3 it follows that $\mathbf{V}_{\mathrm{G}}^* = \mathbf{V}_{\mathrm{NI}}^*$. The result follows by solving (25) for $\mathbf{V}_{\mathrm{NI}}^*$ and substituting the result in (26).

(iv) Lemmas 2(ii) and 3(i) show that $Q(\hat{\theta})$ is a quadratic function of $\sqrt{n}\mathbf{w}$, and the result follows from Lemma 3(ii) and the known result that $Q(\hat{\theta}) \xrightarrow{d} \chi_q^2$ under Assumption 1B.   □

Theorem 1(i) and (iv) actually extend Theorem 1, proved by Satorra [33] for independent groups, to correlated populations and it can be applied to any type of correlated data that can be grouped into a few groups with uncorrelated data (e.g., in panel data by grouping the occasions).

To derive large sample results for $\hat{\theta}$ minimizing (4) under the model in (1) and Assumption 1, we consider the case where all $n^{(i)}$ increase to infinity at a common rate and use $n_m$ as the index for taking a limit in Assumption 2. Assumption 3 is a standard identification condition used in Lemma 1. Note that the true value of $\theta$ in case A of Assumption 1 with fixed variables depends on $n^{(i)}$, since it contains $\bar{\zeta}^{(i)}$ and $\mathbf{S}_{\zeta^{(i)}}$. Thus, we denote the limit of the true value as $\theta_0$. Lemma 1 gives the consistency of the estimator $\hat{\theta}$ that minimizes (4) for the model in (1). Hence, under very weak distributional specifications in Assumption 1, the estimator $\hat{\theta}$ is consistent for the limiting true value $\theta_0$. In fact, it is clear from the proof that the consistency of $\hat{\theta}$ holds for any general mean and covariance structure model $\gamma(\theta) = E(\mathbf{c})$ satisfying $\mathbf{c} \xrightarrow{p} \gamma(\theta_0)$. To characterize the limiting behavior of $\hat{\theta}$ in more detail, especially for the assessment of the so-called asymptotic robustness properties, it is convenient to consider an expansion of $\hat{\theta}$, not around the true value or the limiting true value $\theta_0$, but around some other quantity $\bar{\theta}_n$ defined in Lemma 2 that depends on the unobservable sample moments of the non-normal latent variables $\zeta^{(i)}$ and $\varepsilon_\ell^{(i)}$ ($\ell = 1, \ldots, L^{(i)}$). Thus, the limiting true value $\upsilon_0$ that consists of the true covariance matrices of the random latent variables is replaced in $\bar{\theta}_n$ by $\upsilon_n$ that consists of the unobservable sample moments. While statistical inference is to be made for the true value of $\theta$, $\bar{\theta}_n$ with an artificial quantity $\upsilon_n$ plays a useful role in assessing the property of $\hat{\tau}$ in $\hat{\theta}$, as well as in characterizing the limiting distribution of the whole $\hat{\theta}$ without specifying any moments for $\zeta^{(i)}$ and $\varepsilon_\ell^{(i)}$ ($\ell = 1, \ldots, L^{(i)}$) higher than second order. To obtain an expansion of $\hat{\theta}$ around $\bar{\theta}_n$, we need some smoothness conditions for $\beta^{(i)}(\tau), B^{(i)}(\tau)$ and $\Sigma_{\varepsilon_0^{(i)}}(\tau)$, and the full-column rank of the Jacobian matrix $\mathbf{J}[\gamma(\theta_0)]$ that are stated in Assumption 4. Since the linear independence of the columns of $\mathbf{J}[\gamma(\theta_0)]$ associated with the $\upsilon$ part of $\theta$ is trivial, we need to assume only that the $\tau$ part of the model is specified without redundancy. Thus in Assumption 4 we just assume that $\partial\gamma(\theta_0)/\partial\tau'$ is of full-column rank and Lemma 2 expresses the leading term of $\sqrt{n}(\hat{\theta} - \bar{\theta}_n)$ in terms of $\mathbf{c} - \gamma(\bar{\theta}_n)$. Note that the use of $\bar{\theta}_n$ in Lemma 2 produces an expansion of $\hat{\theta}$ around $\bar{\theta}_n$ with the existence of only second moments of $\zeta^{(i)}$ and $\varepsilon_\ell^{(i)}$ ($\ell = 1, \ldots, L^{(i)}$). It can be shown from the proof that the expansion in Lemma 2 holds for the general model $\gamma(\theta) = E(\mathbf{c})$ and for any $\bar{\theta}_n$ with $\bar{\theta}_n \xrightarrow{p} \theta_0$ provided that $\sqrt{n}[\mathbf{c} - \gamma(\bar{\theta}_n)]$ converges in distribution. However, the special choice of $\bar{\theta}_n$ for the model in (1) makes the result of Lemma 2 practically meaningful. Lemma 3 is actually the key tool in the proof that shows asymptotic robustness. It expresses $\sqrt{n}[\mathbf{c} - \gamma(\bar{\theta}_n)]$ in terms of $\sqrt{n}\mathbf{w}$, which has the same limiting distributions under Assumptions 1 and 1B. Thus, the main difficulty in the proof of Theorem 1(i) is to express $\sqrt{n}(\hat{\tau} - \tau_0)$ in terms of a vector $\sqrt{n}\mathbf{w}$ whose limiting distribution does not depend on the existence of fixed, nonnormal and correlated-over-population variables. Similarly, we proved Theorem 1(iv) by expressing $Q(\hat{\theta})$ as a quadratic function of $\sqrt{n}\mathbf{w}$. The formulas in (5) and (7) in Theorem 1 show what corrections should be made when we have fixed variables in order to get correct asymptotic standard errors for $\hat{\mu}_{\zeta^{(i)}}$ and $\mathrm{vec}(\widehat{\Sigma}_{\zeta^{(i)}})$. These results are novel even for the case with one population. The formula (6) in Theorem 1(ii)(2) shows that the asymptotic standard errors for $\hat{\mu}_{\zeta^{(i)}}$ are robust. Equation (8) in Theorem 1(iii)(2) gives the limiting covariance matrix for $\text{vec}(\widehat{\Sigma}_{\zeta^{(i)}})$ when $\zeta^{(i)}$ are random. Formula (8) involves the computation of fourth-order cumulants of the latent variables $\zeta^{(i)}$ in practice. This is possible in practice and we obtain satisfactory results when we use the errors-in-variables parameterization and have normal errors. For instance, in Example 1 for the model in (2) with normal errors the fourth-order cumulants for $\zeta^{(i)}$ are equal to the fourth-order cumulants of the observed variables for $x^{(i)}$, since the fourth-order cumulants of the normal errors are equal to 0. This technique was used in our simulation study and the results are illustrated in the next section. Note that in most practical cases the measurement errors follow a normal distribution.

Although the paper refers to the multisample case the same theory and methodology can be applied to longitudinal data. That is, two different applications, correlated populations and panel data, can be considered by fitting the same kind of modeling and applying the results presented in this paper. A similar method developed for longitudinal data, known as the general estimating equation (GEE) method, was proposed by Liang and Zeger [19]. The GEE method was proposed for generalized linear models with univariate outcome variables. In this paper several response variables are observed and their relationships are explained by a few latent variables within the time points. It can be shown that a special case of the GEE method, using the identity matrix as the "working" correlation matrix, is a special case of the model considered in this paper. This can be done by treating the outcome variable and the covariates of the generalized linear models as observed variables in the model considered in this paper and setting latent variables equal to covariates by fixing error variances equal to zero. Thus, the results presented in this paper can be also applied to simpler models such as generalized linear models for longitudinal data. On the other hand, the use of a "working" correlation matrix as the one used in the GEE method, could be also used in this methodology in order to increase the efficiency of the method.

Now we define a generalized version of the so-called sandwich estimator used by the GEE method for generalized linear models with the identity matrix as the "working" correlation matrix, and also used by Satorra [28, 29, 30, 31, 32, 33] for latent variable models. We generalize this matrix for correlated populations and we are going to compare it with our proposed matrix $\mathbf{V}_{\text{G}}^{(\theta)}$ defined in Theorem 1 theoretically and numerically. A generalized version of the sandwich (S) estimator is

$$(27) \qquad \mathbf{V}_{\text{S}}^{(\theta)} = \mathbf{A}_0 E(\mathbf{S}_d) \mathbf{A}_0',$$

where $\mathbf{A}_0$ is defined in (11) and $E(\mathbf{S}_d)$ is the expected mean of the sample matrix $\mathbf{S}_d$ that involves third- and fourth-order sample moments defined as

$$\mathbf{S}_d = \begin{pmatrix} \dfrac{1}{n^{(11)}}\mathbf{S}_d^{(11)} & \cdots & \dfrac{1}{n^{(1I)}}\mathbf{S}_d^{(1I)} \\ \vdots & \ddots & \vdots \\ \dfrac{1}{n^{(I1)}}\mathbf{S}_d^{(I1)} & \cdots & \dfrac{1}{n^{(II)}}\mathbf{S}_d^{(II)} \end{pmatrix}$$

with

$$\mathbf{S}_d^{(ik)} = \frac{1}{n^{(ik)}-1}\sum_{j=1}^{n^{(ik)}}(\mathbf{d}_j^{(i)}-\bar{\mathbf{d}}^{(i)})(\mathbf{d}_j^{(k)}-\bar{\mathbf{d}}^{(k)})'$$

and

$$\mathbf{d}_j^{(i)} = \begin{pmatrix} \nu_j^{(i)} \\ \operatorname{vec}[(\nu_j^{(i)}-\bar{\nu}^{(i)})(\nu_j^{(i)}-\bar{\nu}^{(i)})'] \end{pmatrix},$$

where $i,k=1,\ldots,I, j=1,\ldots,n^{(i)}$, and $n^{(ik)}$ denotes the number of correlated individuals between the $i$th and the $k$th populations. Note that the form of the matrix $\mathbf{V}_{\mathrm{S}}^{(\theta)}$ in (27) can be derived from Lemma 2. Equation (12) in Lemma 2 also holds if we replace $\bar{\theta}_n$ by the true value of $\theta$, and the result follows by noting that $\operatorname{Var}[\mathbf{c}-\gamma(\theta_0)] = E(\mathbf{S}_d)$. Theorem 1 actually gives an alternative form of some of the parts of the matrix $\mathbf{V}_{\mathrm{S}}^{(\theta)}$. The parts of the matrix $\mathbf{V}_{\mathrm{G}}^{(\theta)}$ defined in Theorem 1 are actually theoretically exactly the same as the corresponding parts of the matrix $\mathbf{V}_{\mathrm{S}}^{(\theta)}$. In practice, the matrix $\mathbf{A}_0 = \mathbf{A}(\theta_0)$ is estimated by $\hat{\mathbf{A}} = \mathbf{A}(\hat{\theta})$ and the matrix $E(\mathbf{S}_d)$ is estimated by $\mathbf{S}_d$. Despite the fact that the two matrices $\mathbf{V}_{\mathrm{G}}^{(\theta)}$ and $\mathbf{V}_{\mathrm{S}}^{(\theta)}$ are theoretically equal in practice, the asymptotic standard errors (a.s.e.'s) computed by the matrix $\mathbf{V}_{\mathrm{G}}^{(\theta)}$ have less variability than the a.s.e.'s computed by the matrix $\mathbf{V}_{\mathrm{S}}^{(\theta)}$. This happens because the estimation of $\mathbf{V}_{\mathrm{S}}^{(\theta)}$ involves third- and fourth-order moments that are more variable than the second moments of the matrix $\mathbf{V}_{\mathrm{G}}^{(\theta)}$. The matrix $\mathbf{V}_{\mathrm{G}}^{(\theta)}$ involves fourth moments only in the formula of Theorem 1(iii)(2), but these moments do not affect the computation of the other a.s.e.'s. This advantage of using the matrix $\mathbf{V}_{\mathrm{G}}^{(\theta)}$ is shown in the simulation study in the next section.

**4. Simulation study.** We simulate the model in (2) of Example 1. A sample from both populations was generated 1000 times. The simulation was done twice: once with fixed $\zeta^{(i)}$ and once with random $\zeta^{(i)}$ (cases A and B of Assumption 1, respectively). In both cases, $\zeta_j^{(1)}$ and $\zeta_j^{(2)}$ are related (correlated over populations) and were generated as linear combinations of chi-square random variables with 10 degrees of freedom. In case A, a sample of $(\zeta_j^{(1)}, \zeta_j^{(2)})$ was generated with sample means, variances and covariance $\bar{\zeta}^{(1)} = 4.95, \bar{\zeta}^{(2)} = 9.95, s_{\zeta^{(1)}}^2 = 1.97, s_{\zeta^{(2)}}^2 = 1.95$ and $s_{\zeta^{(1)}\zeta^{(2)}} = 1.36$, respectively, and the set of $(\zeta_j^{(1)}, \zeta_j^{(2)})$ was used in all 1000 Monte Carlo samples. In case B, 1000 independent samples were generated for $\{\zeta_j^{(1)}, j = 1, \ldots, 1000; \zeta_j^{(2)}, j = 1, \ldots, 500\}$. The true means, variances and covariance of $\zeta_j^{(1)}$ and $\zeta_j^{(2)}$ are $\mu_{\zeta^{(1)}} = 5, \mu_{\zeta^{(2)}} = 10, \sigma_{\zeta^{(1)}}^2 = 2, \sigma_{\zeta^{(2)}}^2 = 2$ and $\sigma_{\zeta^{(1)}\zeta^{(2)}} = 1.4$. Note that the above means and variances are estimated, but the covariance $\sigma_{\zeta^{(1)}\zeta^{(2)}}$ is not, in accordance with the estimation method that we suggest. Note that we suggest this method for several populations with quite unbalanced data. In this study it is easy to use the full likelihood and estimate the covariance $\sigma_{\zeta^{(1)}\zeta^{(2)}}$, but this is not always true in more complicated studies. By not estimating some of the covariances between the two populations, we lose some efficiency, for example, we obtain larger a.s.e.'s. We discuss the efficiency of the method in more detail later in this section.

In both cases A and B, 1000 samples were generated for independent $e_\ell^{(i)}, i = 1, 2, \ell = 0, 1, \ldots, L^{(i)}$, with $L^{(1)} = 3$ and $L^{(2)} = 2$. The errors $e_{0j}^{(i)}, i = 1, 2$, are normally distributed with mean 0 and unknown variance $\sigma_{e_0^{(i)}}^2$, while all the other errors $e_{\ell j}^{(i)}$ for $i = 1, 2, \ell = 1, \ldots, L^{(i)}$, were generated from a chi-square distribution with 10 degrees of freedom, $\chi_{10}^2$, with adjusted mean 0 and variance $\sigma_{e_{0j}^{(i)}}^2$. The variance for $e_{0j}^{(i)}$ is common for the two populations, $\sigma_{e_0}^2 = \sigma_{e_0^{(1)}}^2 = \sigma_{e_0^{(2)}}^2$. In both cases with fixed and random $\zeta_j^{(i)}$, the true values for the error variances are $\sigma_{e_0^{(i)}}^{o2} = \sigma_{e_1^{(i)}}^{o2} = 0.1$ and $\sigma_{e_2^{(i)}}^{o2} = \sigma_{e_3}^{o2} = 0.2$, and the true value for the vector $\tau$ is $\tau^0 = (1, 2, -1, -0.1, 0.1, -0.01, 1, 0.1)$. The parameters of $\tau$ are shown in the first column of the first part of Table 1. In accordance with the notation of this paper, the vector $\theta = (\tau', \upsilon')'$, where $\upsilon$ contains $\sigma_{e_\ell^{(i)}}^2$ $(i = 1, 2, \ell = 1, \ldots, L^{(i)})$ and the means and variances of $\zeta_j^{(i)}$ $(i = 1, 2)$. To estimate $\theta$, we use normal MLE by minimizing (4) despite the appearance of fixed and nonnormal variables, and when we estimate the parameters, we are pretending that we do not know the true values of the parameters.

Some of the results in the simulation study are shown in the first part of Table 1. Columns 2, 4 and 6 show results from case A with fixed $\zeta_j^{(i)}$, while columns 3, 5 and 7 show results from case B with random $\zeta_j^{(i)}$. Columns 2 and 3 of Table 1 compare the a.s.e.'s (Gse) computed by the matrix $\mathbf{V}_G^{(\tau)}$ in

TABLE 1
*Results from the sumulation study**

| | | Bias of Gse $\frac{\text{Gse}}{\text{MCse}}$ | | Variability of Gse $\frac{\text{SMCse}}{\text{GMCse}}$ | | Efficiency of the method relative to the full likelihood $\frac{\text{PL}-\text{MCse}}{\text{FL}-\text{MCse}}$ | |
|---|---|---|---|---|---|---|---|
| Parameters $\tau$ | | Fixed | Random | Fixed | Random | Fixed | Random |
| $\beta_1$ | | 1.01 | 1.01 | 1.63 | 1.56 | 0.99 | 1.03 |
| $\beta_2$ | | 1.01 | 0.99 | 1.78 | 1.68 | 1.01 | 1.05 |
| $\beta_3$ | | 0.97 | 1.00 | 1.84 | 1.50 | 1.00 | 1.06 |
| $\gamma_1$ | | 1.00 | 0.99 | 1.44 | 1.47 | 1.00 | 1.04 |
| $\gamma_2$ | | 0.97 | 0.99 | 2.02 | 1.56 | 1.01 | 1.05 |
| $\delta_1$ | | 1.00 | 1.00 | 1.65 | 1.57 | 1.00 | 1.03 |
| $\delta_2$ | | 1.00 | 0.98 | 1.60 | 1.44 | 1.02 | 1.06 |
| $\sigma_{e_0}^2$ | | 0.99 | 0.99 | 2.68 | 1.56 | 1.00 | 1.03 |
| Results for $\gamma_1$ under different distribution assumptions—degrees of freedom for $\zeta_j^{(i)} \sim \chi^2(d_1)$ and $e_{k,j}^{(i)} \sim \chi^2(d_2)$ | | | | | | | |
| $d_1$ | $d_2$ | | | | | | |
| 1 | 1 | 1.00 | 1.00 | 1.59 | 1.69 | 1.01 | 1.09 |
| 3 | 3 | 1.00 | 1.01 | 1.55 | 1.43 | 1.01 | 1.07 |
| 3 | 10 | 0.99 | 0.98 | 1.48 | 1.41 | 1.01 | 1.07 |
| 10 | 3 | 0.99 | 1.00 | 1.51 | 1.51 | 1.01 | 1.04 |

* Monte Carlo standard errors (MCse) for the estimated parameters in $\tau$ versus the proposed a.s.e.'s (Gse) of $\hat{\tau}$, computed by $\mathbf{V}_{\mathrm{G}}^{(\tau)}$ defined in Theorem 1. Comparison between the MCse for Gse (GMCse) and the MCse for the a.s.e.'s computed by the sandwich estimator, $\mathbf{V}_{\mathrm{S}}^{(\theta)}$, given in (27) (SMCse). MCse computed under the full likelihood (FL) and under the partial likelihood (PL). Results are shown for cases A and B of Assumption 1 with fixed and random $\zeta_j^{(i)}$.

Theorem 1(i) with the Monte Carlo standard errors (MCse). All the ratios are 1 or very close to 1 and this means that the proposed a.s.e.'s have very small bias. Bias exists because we use the a.s.e.'s as estimates for the true s.e.'s of the parameters in finite samples. Actually, Lemma 1 proves that the bias converges to 0 as the sample sizes increase to infinity. In this study, for sample sizes $n^{(1)} = 1000$ and $n^{(2)} = 500$, the bias is negligible.

Now we compute Monte Carlo standard errors for the a.s.e.'s computed by the matrix $\mathbf{V}_{\mathrm{G}}^{(\theta)}$ (GMCse) and for the a.s.e.'s computed by the matrix $\mathbf{V}_{\mathrm{S}}^{(\theta)}$ (SMCse), defined in (27). The ratio (SMCse)/(GMCse) compares the variability of the two different estimates of the a.s.e.'s. This ratio is computed for the parameters in $\tau$ and the results are shown in columns 4 and 5 of Table 1 for both cases with fixed and random $\zeta_j^{(i)}$. All the ratios are significantly larger than 1 and this fact indicates that the a.s.e.'s computed by the sandwich estimator $\mathbf{V}_S^{(\theta)}$ have larger variability than the a.s.e.'s computed by our suggested estimator $\mathbf{V}_G^{(\theta)}$.

Now, as to the efficiency of the method, we computed the a.s.e.'s under the full likelihood (FL) and under the partial likelihood (PL) given in (4). The ratio of the two a.s.e.'s,

$$\text{(28)} \qquad \text{efficiency} = \frac{\text{PL} - \text{MCse}}{\text{FL} - \text{MCse}},$$

is given for all the parameters in $\tau$ in the last two columns of Table 1. These ratios actually show the efficiency of the method relative to the FL. In both cases the efficiency is very satisfactory since the ratios are close to 1. The efficiency loss is very small for case A with fixed $\zeta_j^{(i)}$ and relatively small for case B with random $\zeta_j^{(i)}$.

In the second part of Table 1, we make the nonnormal distributions more skewed to the right by changing the degrees of freedom, $d_1$, and $d_2$, for $\zeta_j^{(i)} \sim \chi^2(d_1)$ and $e_{k,j}^{(i)} \sim \chi^2(d_2)$. All the results remain the same for case A with fixed $\zeta_j^{(i)}$ and the only difference in case B with random $\zeta_j^{(i)}$ is that the efficiency ratio of the method relative to the full likelihood (last column) becomes larger but remains smaller than 1.10 even in the extreme case with 1 degree of freedom for both $d_1$ and $d_2$. Thus, the derived asymptotic standard errors give satisfactory results for distributions with very long tails that often appear in applications (e.g., in finance and banking).

For the parameters $\mu_{\zeta^{(1)}}, \mu_{\zeta^{(2)}}, \sigma_{\zeta^{(1)}}^2$ and $\sigma_{\zeta^{(2)}}^2$ we used the formulas in (5), (6), (7) and (8) provided in Theorem 1(ii) and (iii) and we derived results similar to the previous ones. It should be pointed out that the sandwich estimator does not provide correct a.s.e.'s for case A with fixed $\zeta_j^{(i)}$ for the parameters $\mu_{\zeta^{(1)}}, \mu_{\zeta^{(2)}}, \sigma_{\zeta^{(1)}}^2$ and $\sigma_{\zeta^{(2)}}^2$. Our novel formulas in (5) and (7) show what corrections should be made in order to obtain the correct a.s.e.'s in this case. The a.s.e.'s are evaluated at the estimated value of $\theta$, $\hat{\theta}$. Note that all the a.s.e.'s are functions of $\theta$ except the ones for $\hat{\sigma}_{\zeta^{(1)}}^2$ and $\hat{\sigma}_{\zeta^{(2)}}^2$ (elements of the matrix $\widehat{\Sigma}_{\zeta^{(i)}}$ in Theorem 1) that require fourth moments (or cumulants) for $\zeta_j^{(i)}$. In general, the fourth-order cumulants, $\psi$, are prescribed by the following property: if $x = y + z$ with $y$ and $z$ independent random variables, then $\psi_x = \psi_y + \psi_z$. Thus, in the model used in the simulation, it holds that $\psi_{x^{(i)}} = \psi_{\zeta^{(i)}} + 0$, since the errors, $e_{0j}^{(i)}$, are assumed to be normal, having fourth-order cumulants equal to 0. Thus, the sample fourth-order cumulants of $x^{(i)}$ were used for the computation of the a.s.e.'s for $\hat{\sigma}_{\zeta^{(1)}}^2$ and $\hat{\sigma}_{\zeta^{(2)}}^2$.

The a.s.e.'s can be used for hypothesis testing of the parameters. The power of the tests is also robust when the sample sizes are quite large due to the applicability of the multivariate central limit theorem. In the above simulation study, we use, as an example, $H_0 : \delta_1 = 0$ versus $H_1 : \delta_1 < 0$ in case A with fixed $\zeta_j^{(i)}$. Using level of significance $\alpha = 0.05$, $H_0$ is rejected when $z < -1.645$ where $z = \hat{\delta}_1 / \hat{\sigma}_{\delta_1}^2$. Thus, the expected power (EP) is approximately

$$(29) \qquad \text{EP}(\delta_1^*) = \Phi\left(-1.645 + \frac{\delta_1^*}{\text{MCse for } \delta_1}\right) = 0.956,$$

where the $\Phi$ function is the standard cumulative normal distribution and we compute the power for the actual value of $\delta_1$, $\delta_1^* = -0.01$. We also compute the simulated power (SP) as

$$(30) \qquad \text{SP} = \frac{\# \text{ of times that } [\hat{\delta}_1 / (\text{a.s.e. of } \hat{\delta}_1)] < -1.645}{1000} = 0.967.$$

Thus, the results support the robustness of power for nonnormal and correlated populations. The power for overall-fit measures was investigated by Satorra and Saris [36] and Satora [34] in structural equation models.

The robustness of the chi-square test statistic is shown in Table 2 for case A with fixed $\zeta_j^{(i)}$. The mean and the variance of the 1000 simulated values of $Q(\hat{\theta})$ in (4) are close to the expected 6 and 12, respectively. Also, the simulated percentiles in the second row are close to the expected ones given in the first row of Table 2. For similar studies using simpler models, see [30, 32, 33] and [25].

In summary, the model in (1) with the errors-in-variables parameterization can formulate the multipopulation analysis in a meaningful fashion. The corresponding statistical analysis under the pseudo-normal-independence model gives a simple and correct way to conduct statistical inferences about the parameter vector $\tau$ without specifying a distributional form or dependency structure over populations. In practice, $\tau$ contains all the parameters of direct interest. The asymptotic covariance matrix and standard errors can be readily computed using existing procedures, and provide a good approximation in moderately sized samples. The proposed a.s.e.'s have smaller variability than the variability of the robust sandwich estimator, provide high efficiency relative to the full-likelihood method and can be used for hypothesis testing with robust power. For instance, in the simulation study for one of the most important parameters, $\delta_1$, in case A with fixed $\zeta_j^{(i)}$, the variability ratio is 1.65 (see Table 1), the efficiency ratio is 1.00 (see Table 1) and the power of the test $H_0 : \delta_1 = 0$ versus $H_1 : \delta_1 < 0$ is 0.967. That is, if the standard deviation of our proposed a.s.e. for $\delta_1$ is 1, then the standard deviation of $\delta_1$ computed by the robust sandwich estimator is 1.65. Also, our proposed a.s.e. for $\delta_1$ is close enough to the a.s.e. for $\delta_1$ using the full likelihood, and the power of the test is very high, 0.967, and very close to the expected power, 0.960.

Table 2
*Monte Carlo mean, variance and percentiles for the chi-square test statistics with 6 degrees of freedom*

| Mean = 6 | Variance = 12 | 10% | 25% | 50% | 75% | 90% | 95% | 99% |
|---|---|---|---|---|---|---|---|---|
| 6.0 | 11.7 | 9.2 | 23.6 | 49.7 | 75.9 | 90.5 | 96.3 | 98.9 |

**5. Application.** An application for model (1), estimated by minimizing (4), and for Theorem 1 is presented by analyzing a data set from the Bank of Greece with annual statements for the period 1999–2003. We examine the relationship between asset risk and capital in the Greek banking sector. As capital, we use total capital over total bank assets (capital-to-asset ratio). The variable for total capital includes core capital (tier I) plus supplementary capital (tier II) minus deductions such as holdings of capital of other credit and financial institutions. As measures for asset risk, we use the two main components of risk-weighted assets which reflect credit and market risk. There is a two-way direction effect between capital and asset risk, and these relationships can be analyzed in a multivariate setting with simultaneous equations; see [7] for the life insurance industry. Baranoff, Papadopoulos and Sager [6] compared the effect of two measures for asset risk to capital structure by approaching latent variables for the risk factors via a dynamic structural equation model, and they pointed out the differences between large and small companies. They fitted latent variable models on a balanced data set concentrating on companies for which data for all years are available. In such balanced cases we ignore companies that have been bankrupt or have been merged with other companies, and new companies that started after the first year. In many studies, researchers are interested in examining such companies and fit latent variables, such as macroeconomic and risk factors or measurement errors, in a highly unbalanced data set. Papadopoulos and Amemiya [26] discussed the disadvantages of the existing methods for unbalanced data. The methodology proposed in this paper is appropriate for highly correlated, nonnormal and unbalanced data. Also, Theorem 1 ensures robust asymptotic standard errors and overall-fit measures.

In this paper we analyze first differences of the logarithmic (ln) transformation, which actually approximate percentage changes, in order to avoid spurious regression, nonstationarity and cointegration to some extent. The explicit form of the model is

$$\Delta \ln\left(\frac{\text{capital}}{\text{assets}}\right)_j^{(t)} = \beta_1 \zeta_j^{(t)} + \varepsilon_{1j}^{(t)},$$

$$\Delta \ln \left( \frac{\text{credit risk}}{\text{assets}} \right)^{(t)}_j = \beta_2 \zeta^{(t)}_j + \varepsilon^{(t)}_{2j},$$

$$(31) \qquad t = 2000, \ldots, 2003, j = 1, 2, \ldots, n^{(t)},$$

$$\Delta \ln \left( \frac{\text{market risk}}{\text{assets}} \right)^{(t)}_j = \beta_3 \zeta^{(t)}_j + \varepsilon^{(t)}_{3j}.$$

The above model is a confirmatory factor analytic model with one underlying factor, $\zeta^{(t)}_j$, that explains the relationships of the three observed variables, and it is a simple case of model (1). The parameter $\beta_1$ is fixed equal to 1, for identification reasons, and this actually assigns the latent factor, $\zeta^{(t)}_j$, to have the same units as the corresponding observed variable. The variables $\zeta^{(t)}_j, \varepsilon^{(t)}_{2j}$ and $\varepsilon^{(t)}_{3j}$ are assumed to follow nonnormal distributions, since the observed variables have long tails, which is very common for financial variables. These variables also have unrestricted variances over time due to the heteroskedasticity over time of the observed variables. By viewing $\varepsilon^{(t)}_{1j}$ as measurement error, then as a smooth and invariant latent variable over time it is assumed to follow a normal distribution with equal variances over time. Also, we assume that the autocorrelation of the observed variables is explained by the autocorrelation of $\zeta^{(t)}_j$ and that the errors $\varepsilon^{(t)}_{kj}, k = 1, 2, 3$, are independent over time, which is a common assumption when we analyze differences and applications in this analysis. In general, if there is still autocorrelation after taking the first differences, we can try second differences, and so on.

Frequently, in finance and banking we are interested in examining the relationship between asset risk and capital ratio, particularly when the asset risk increases or decreases significantly. In these cases the restricted variables of asset risk have truncated distributions, in addition to their long tails, and the issue of robustness of standard methods to such nonnormal data becomes very important and necessary. Especially in the cases with restricted variables, the already unbalanced data lose the appearance of the banks in consecutive years, since they do not satisfy the required condition every year. Therefore, it is very difficult and in many, if not all, applications it is impossible to model the time series structure. Then methodologies that focus on modeling relationships between variables within the occasions, such as the proposed model in (1), become very attractive and useful.

Table 3 shows results for model (31) using the proposed methodology for all data and for data arising by restricting one of the observed variables. For more details, see the explanation in Table 3. Table 4 shows the explicit pattern of missing values for the case with market asset risk less than $-0.05$. Thus, if we try to reformulate the four correlated samples as independent samples based on the missing pattern of the banks, then we end up with

TABLE 3

*Results for the coefficients $\beta_k$, $k = 1, 2, 3$, of model (31) for several cases: for all available data (column 2) and for data that arise by restricting one of the observed variables to be significantly positive ($> 0.05$) (columns 3, 5 and 7) or be negative ($< -0.05$) (columns 4, 6 and 8)*[*]

|  | **Without restrictions** | **Restrictions on capital-to-asset ratio** | | **Restrictions on credit risk ratio** | | **Restrictions on market risk ratio** | |
|---|---|---|---|---|---|---|---|
|  | **All** $n = 68$ | **> 0.05** $n = 23$ | **< −0.05** $n = 39$ | **> 0.05** $n = 37$ | **< −0.05** $n = 18$ | **> 0.05** $n = 26$ | **< −0.05** $n = 41$ |
| $\beta_1$ | 0.96 | 0.53 | 0.47 | 0.61 | 1.00 | 0.82 | 1.00 |
|  | 1.00 | 1.00 | 1.00 | 1.00 | 1.00 | 1.00 | 1.00 |
|  | (—) | (—) | (—) | (—) | (—) | (—) | (—) |
| $\beta_2$ | 0.45 | 0.36 | 0.57 | 0.16 | −0.03 | 0.54 | 0.58 |
|  | 0.46 | 0.68 | 1.21 | 0.25 | −0.03 | 0.66 | 0.58 |
|  | (1.95) | (1.58) | (0.43) | (0.94) | (−0.13) | (2.18) | (4.57) |
| $\beta_3$ | 0.48 | 1.00 | 0.16 | 1.00 | 0.54 | 0.73 | 0.51 |
|  | 0.50 | 1.88 | 0.34 | 1.64 | 0.54 | 0.89 | 0.51 |
|  | (1.98) | (3.00) | (0.56) | (4.69) | (2.74) | (2.37) | (3.76) |

[*]For each cell we report the standardized (first row; see [10] for a definition) and the unstandardized (second row) coefficients, and the value of the $z$ test [unstandardized coefficient over its asymptotic standard error (a.s.e.)]. The sum of the sample sizes for the four years, $n^{(2000)} + n^{(2001)} + n^{(2002)} + n^{(2003)}$, appears in the third row for each case.

TABLE 4

*Pattern of missing data for the case with differences of the* ln*'s for market risk ratio* $< -0.05$[*]

| Group | Number of banks | 2000 | 2001 | 2002 | 2003 |
|---|---|---|---|---|---|
| 1 | 2 | 0 | 0 | 0 | 2 |
| 2 | 1 | 0 | 0 | 1 | 0 |
| 3 | 2 | 0 | 0 | 2 | 2 |
| 4 | 4 | 0 | 4 | 0 | 4 |
| 5 | 1 | 0 | 1 | 1 | 0 |
| 6 | 3 | 0 | 3 | 3 | 3 |
| 7 | 1 | 1 | 0 | 1 | 0 |
| 8 | 1 | 1 | 0 | 1 | 1 |
| 9 | 1 | 1 | 1 | 0 | 1 |
| 10 | 1 | 1 | 1 | 1 | 0 |
| 11 | 1 | 1 | 1 | 1 | 1 |
| Total number of banks | 18 | 5 | 11 | 11 | 14 |

[*]In the last four columns the nonzero numbers indicate that for the corresponding group (numbered in column 1) the number of banks stated in column 2 appears in these particular years labelled in row 1. The nonzero numbers in columns 2–6 are the same in each row.

11 independent samples that have very small sample sizes—smaller than four—and most of them having just one observation. Therefore, the analysis of balanced data is not possible since there is only one bank that appears in all four years that satisfies the required restriction. Also, the analysis of time series structure is not possible, since all samples that have banks appearing in any two or more consecutive years have sample sizes less than three. The methodology suggested in this paper can be applied to four correlated samples with observations from the four years, respectively. The sample sizes for the four years are 5, 11, 11 and 14 from 2000, 2001, 2002 and 2003, respectively, and the sum of the four samples is 41 (see the last row in Table 4). According to our methodology, we analyze 41 observations from banks that appear in at least one year. In this case, there are 18 different banks that appear in some of the four years. It should be noted that the estimated parameters of interest, $\beta_2$ and $\beta_3$, belong to the vector $\tau$ and thus, according to Theorem 1(i), their asymptotic standard errors can be computed by the covariance matrix $\mathbf{V}_{\mathrm{NI}}^{(\tau)}$. The computation of $\mathbf{V}_{\mathrm{NI}}^{(\tau)}$ involves moments only of first and second order, and this issue is very important especially when the sample size, as in this example, is small. Only the asymptotic covariance matrix $\mathbf{V}_{\mathrm{G}}^{(\mathrm{vec}(\Sigma_{\zeta^{(i)}}))}$, defined in (8), requires fourth-order moments for its computation, and for its use we need larger sampler sizes than the sample sizes of this example. Thus, we can fit panel data models of moderate sample sizes relative to the number of estimated parameters and make statistical inference for the most important parameters without using moments of order higher than two in the analysis.

Also, in the case with all banks (with no restriction on any observed variable), there are 20 different banks that provide data for some of the four years, creating a very unbalanced data set with only 12 banks appearing in all four years. As Table 3 shows in this case, if we add the banks that appear every year, then we have a total of 68 observations from the 20 banks. Actually, these 68 observations were analyzed in four correlated samples, giving consistent estimates, and correct and efficient asymptotic standard errors relative to the sandwich estimator, despite the nonnormality and autocorrelation of the variables, according to Theorem 1.

The standardized coefficients in Table 3, in the case without restrictions on the observed variables (column 2), indicate that the latent factor, $\zeta_j^{(t)}$, is strongly associated with the capital-to-asset ratio, 0.96, and has almost the same degree of correlation with the two measures for asset risk, 0.45 and 0.48. The results significantly change when we restrict one of the observed values on significantly positive or negative. When we restrict the capital-to-asset ratio on positive values, the factor $\zeta_j^{(t)}$ coincides with market risk, and gives a stronger and significant correlation with capital-to-asset ratio than the one with credit risk. The results found by restricting capital-to-asset ratio on negative values are not statistically significant. When we restrict the credit risk ratio on positive and on negative values, the factor $\zeta_j^{(t)}$ coincides with market risk and capital-to-asset ratio, respectively, and is significantly correlated with capital-to-asset ratio and market risk, respectively, 0.61 and 0.54, and not with the other variable. Comparing the results from the last two columns to the results of column 2, we observe that the standardized coefficients for $\beta_2$ and $\beta_3$ are higher in these columns than the ones in column 2. Also note that in column 7 the market risk gives a much higher standardized coefficient, 0.73, than the credit risk, 0.54, and indicates the strongest relationship between capital-to-asset ratio and asset risk. All in all, as expected, the capital-to-asset ratio is always positively correlated to both credit and market asset risk. Also, the results change when we restrict one of the observed variables to be positive or negative, and thus it is worthwhile. Even though the panel data are highly unbalanced and additionally lose their consecutive appearance over the years, our methodology can be applied and can provide correct statistical inference.

**Acknowledgments.** The authors wish to thank the reviewers and editors for their insightful and constructive comments. The first author especially thanks Director Panagiotis Kiriakopoulos, Professor David Scott and Professor Wayne Fuller.

## REFERENCES

[1] AMEMIYA, Y. and ANDERSON, T. W. (1990). Asymptotic chi-square tests for a large class of factor analysis models. *Ann. Statist.* **18** 1453–1463. MR1062719

[2] AMEMIYA, Y., FULLER, W. A. and PANTULA, S. G. (1987). The asymptotic distributions of some estimators for a factor analysis model. *J. Multivariate Anal.* **22** 51–64. MR0890881

[3] ANDERSON, T. W. (1987). Multivariate linear relations. In *Proc. Second International Tampere Conference in Statistics* (T. Pukkila and S. Puntanen, eds.) 9–36. Univ. Tampere, Finland.

[4] ANDERSON, T. W. (1989). Linear latent variable models and covariance structures. *J. Econometrics* **41** 91–119. MR1007726

[5] ANDERSON, T. W. and AMEMIYA, Y. (1988). The asymptotic normal distribution of estimators in factor analysis under general conditions. *Ann. Statist.* **16** 759–771. MR0947576

[6] BARANOFF, E. G., PAPADOPOULOS, S. and SAGER, T. W. (2005). The effect of regulatory versus market asset risk on the capital structure of life insurers: A structural equation modeling approach. Working paper. Available at smealsearch2.psu.edu/121.html.

[7] BARANOFF, E. G. and SAGER, T. W. (2002). The relations among asset risk, product risk, and capital in the life insurance industry. *J. Banking and Finance* **26** 1181–1197.

[8] BENTLER, P. M. (1983). Some contributions to efficient statistics in structural models: Specification and estimation of moment structures. *Psychometrika* **48** 493–517. MR0731206

[9] BENTLER, P. M. (1989). *EQS Structural Equations Program Manual.* BMDP Statistical Software, Los Angeles.

[10] BOLLEN, K. A. (1989). *Structural Equations with Latent Variables.* Wiley, New York. MR0996025

[11] BROWNE, M. W. (1984). Asymptotically distribution-free methods for the analysis of covariance structures. *British J. Math. Statist. Psych.* **37** 62–83. MR0783499

[12] BROWNE, M. W. (1987). Robustness in statistical inference in factor analysis and related models. *Biometrika* **74** 375–384. MR0903138

[13] BROWNE, M. W. (1990). Asymptotic robustness of normal theory methods for the analysis of latent curves. In *Statistical Analysis of Measurement Error Models and Applications* (P. J. Brown and W. A. Fuller, eds.) 211–225. MR1087111

[14] BROWNE, M. W. and SHAPIRO, A. (1988). Robustness of normal theory methods in the analysis of linear latent variate models. *British J. Math. Statist. Psych.* **41** 193–208. MR0985133

[15] CHOU, C.-P., BENTLER, P. M. and SATORRA, A. (1991). Scaled test statistics and robust standard errors for nonnormal data in covariance structure analysis: A Monte Carlo study. *British J. Math. Statist. Psych.* **44** 347–358.

[16] FULLER, W. A. (1987). *Measurement Error Models.* Wiley, New York. MR0898653

[17] JÖRESKOG, K. (1971). Simultaneous factor analysis in several populations. *Psychometrika* **36** 409–426.

[18] JÖRESKOG, K. and SÖRBOM, D. (1989). *LISREL* 7: *A Guide to the Program and Applications,* 2nd ed. SPSS, Chicago.

[19] LIANG, K. Y. and ZEGER, S. L. (1986). Longitudinal data analysis using generalized linear models. *Biometrika* **73** 13–22. MR0836430

[20] LEE, S. Y. and TSUI, K. L. (1982). Covariance structure analysis in several populations. *Psychometrika* **47** 297–308. MR0678064

[21] MAGNUS, J. and NEUDECKER, H. (1988). *Matrix Differential Calculus with Applications in Statistics and Econometrics.* Wiley, New York. MR0940471

[22] MOOIJAART, A. and BENTLER, P. M. (1991). Robustness of normal theory statistics in structural equation models. *Statist. Neerlandica* **45** 159–171. MR1129198

[23] MUTHÉN, B. (1989). Multiple group structural modeling with nonnormal continuous variables. *British J. Math. Statist. Psych.* **42** 55–61.

[24] MUTHÉN, B. and KAPLAN, D. (1992). A comparison of some methodologies for the factor analysis of nonnormal Likert variables: A note on the size of the model. *British J. Math. Statist. Psych.* **45** 19–30.

[25] PAPADOPOULOS, S. and AMEMIYA, Y. (1994). Asymptotic robustness for the structural equation analysis of several populations. *Proc. Business and Economics Statistics Section* 65–70. Amer. Statist. Assoc., Alexandria, VA.

[26] PAPADOPOULOS, S. and AMEMIYA, Y. (1995). On factor analysis of longitudinal data. *Proc. Biometric Statistics Section* 130–135. Amer. Statist. Assoc., Alexandria, VA.

[27] SAS INSTITUTE, INC. (1990). *SAS/STAT User's Guide, Version 6* **1**, 4th ed. SAS Institute, Cary, NC.

[28] SATORRA, A. (1992). Asymptotic robust inferences in the analysis of mean and covariance structures. *Sociological Methodology* **22** 249–278.

[29] SATORRA, A. (1993). Multi-sample analysis of moment-structures: Asymptotic validity of inferences based on second-order moments. In *Statistical Modeling and Latent Variables* (K. Haagen, D. J. Bartholomew and M. Deistler, eds.) 283–298. North-Holland, Amsterdam. MR1236719

[30] Satorra, A. (1993). Asymptotic robust inferences in multi-sample analysis of augmented-moment structures. In *Multivariate Analysis*: *Future Directions* **2** (C. M. Cuadras and C. R. Rao, eds.) 211–229. North-Holland, Amsterdam. MR1268430

[31] Satorra, A. (1994). On asymptotic robustness in multiple-group analysis of multivariate relations. Paper presented at the Conference on Latent Variable Modeling with Applications to Causality, Los Angeles.

[32] Satorra, A. (1997). Fusion of data sets in multivariate linear regression with errors-in-variables. In *Classification and Knowledge Organization* (R. Klar and O. Opitz, eds.) 195–207. Springer, London.

[33] Satorra, A. (2002). Asymptotic robustness in multiple group linear-latent variable models. *Econometric Theory* **18** 297–312. MR1891826

[34] Satorra, A. (2003). Power of chi-square goodness-of-fit test in structural equation models: The case of nonnormal data. In *New Developments of Psychometrics* (H. Yanai, A. Okada, K. Shigemasu, Y. Kano and J. J. Meulman, eds.) 57–68. Springer, Tokyo.

[35] Satorra, A. and Bentler, P. M. (1990). Model conditions for asymptotic robustness in the analysis of linear relations. *Comput. Statist. Data Anal.* **10** 235–249. MR1086038

[36] Satorra, A. and Saris, W. E. (1985). Power of the likelihood ratio test in covariance structure analysis. *Psychometrika* **50** 83–90. MR0789217

[37] Shapiro, A. (1987). Robustness properties of the MDF analysis of moment structures. *South African Statist. J.* **21** 39–62. MR0903793

Department for the Supervision
 of Credit System and FI
Division of Statistics
Bank of Greece
Amerikis 3
Athens 10250
Greece
e-mail: sapapa@bankofgreece.gr

IBM T. J. Watson Research Center
Route 174
Yorktown Heights, New York 10598
e-mail: yasuo@us.ibm.com